\documentclass{article}
\usepackage{amsmath,amssymb,amsthm,amscd,latexsym,}
 
\usepackage{graphicx}

\makeatletter
\renewcommand{\mod}[1]{\allowbreak \if@display \mkern 8mu \else
\mkern 5mu\fi {\operator@font mod}\,\,#1}
\makeatother

\newcommand{\bn}{\mathbb N}
 \newcommand{\bq}{\mathbb Q}
\newcommand{\br}{\mathbb R}
\newcommand{\bz}{\mathbb Z}
\newcommand{\bff}{\mathbb F}

\newcommand{\bk}{\mathbb K}

\newtheorem{theorem}{Theorem}[section]

\newtheorem{definition}[theorem]{Definition}

\newcommand\F{\mathcal F}
\newcommand\Hh{\mathcal H}

\newcommand\M{\mathcal M}

\begin{document}
\title{On ground fields of arithmetic hyperbolic reflection groups. III}
\date{}
\author{Viacheslav V. Nikulin\footnote{This paper was written with the
financial support of EPSRC, United Kingdom (grant no. EP/D061997/1)}}

\maketitle

\begin{abstract}
This paper continues \cite{Nik5} (arXiv.org:math.AG/0609256),
\cite{Nik6} \linebreak
(arXiv:0708.3991) and \cite{Nik7}(arXiv:0710.0162)

Using authors's methods of 1980, 1981,
some explicit finite sets of number fields containing all ground fields of
arithmetic hyperbolic reflection groups in dimension at least 3
are defined, and explicit bounds of
their degrees (over $\bq$) are obtained.

Thus, now, explicit bound of degree of ground fields of arithmetic
hyperbolic reflection groups is known in all dimensions. Thus, now,
we can, in principle, obtain effective finite classification
of arithmetic hyperbolic reflection groups in all dimensions
together.
\end{abstract}

\vskip1cm \centerline{\it To 85th Birthday of Igor Rostislavich
Shafarevich} \vskip1cm

\section{Introduction} \label{introduction} This paper continues 
\cite{Nik5}, \cite{Nik6} and \cite{Nik7}. 
See introductions of these papers about
history, definitions and results concerning the subject.

In \cite{Nik6} and \cite{Nik7} some explicit finite sets of totally real
algebraic number fields containing all ground fields of arithmetic hyperbolic
reflection groups in dimensions $n\ge 4$ were defined,
and good explicit bounds of degrees (over $\bq$) of their fields were
obtained.  In particular, an explicit bound ($\le 56$) for $n\ge 6$ and
($\le 138$) for $n\ge 4$ of degree of the ground field of any arithmetic
hyperbolic reflection group in dimension $n\ge 4$ was obtained.

Here we continue this study for dimensions $n=3$. Using similar
methods, we define some explicit finite sets of totally real
algebraic number fields containing all ground fields of arithmetic
hyperbolic reflection groups in dimensions $n \ge 3$. Moreover, an
explicit bound ($\le 909$) of degrees of fields from these sets
are obtained. Here we also use result by Long, Maclachlan and Reid
from \cite{LMR}.

Thus, degree of the ground field of any arithmetic hyperbolic
reflection group of dimension $n\ge 3$ is bounded by $909$. Remark
that finiteness of the number of maximal arithmetic hyperbolic
reflection groups of dimension 3 was obtained by Agol \cite{Agol},
it follows a theoretical existence of some bound of the degree of
ground fields of arithmetic hyperbolic reflection groups in
dimension 3. Difference of our result here is that we give an
explicit bound for the degree which is important for finite
effective classification. It also gives another proof
of finiteness in dimension 3. In fact, using our methods, we show
that finiteness in dimension 3 follows from finiteness in
dimension 2.

It is also very important that all these fields are attached to fundamental
chambers of arithmetic hyperbolic reflection groups, and they can be further
geometrically investigated and restricted.

It was also shown in \cite{Nik6} (using results of \cite{LMR} and
\cite{Borel}, \cite{Tak4}) that degree of the ground field of
any arithmetic hyperbolic reflection group of dimension $n=2$ is
bounded by $44$.

Thus, now an explicit bound of degree of ground fields of
arithmetic hyperbolic reflection groups is known in all
dimensions. By \cite{Nik1} and \cite{Nik2} and \cite{Vin2},
\cite{Vin3}, then there exists an effective finite classification
of maximal arithmetic hyperbolic reflection groups in all
dimensions together. More generally, there exists an effective
finite classification of similarity classes of reflective
hyperbolic lattices $S$. Their full reflection groups $W(S)$
contain all maximal arithmetic hyperbolic reflection groups as
subgroups of finite index.

Since this paper is a direct continuation of \cite{Nik6} and
\cite{Nik7}, we use notations, definitions and results of these
papers.


\section{Ground fields of arithmetic hyperbolic reflection
groups in dimensions $n\ge 3$}
\label{sec:5and4}

Since this paper is a direct continuation of \cite{Nik6} and
\cite{Nik7}, we use notations, definitions and results of this
papers.

In \cite[Secs 3 and 4]{Nik6}, explicit finite sets
$\F L^4$, $\F T$, $\F\Gamma^{(4)}_i(14)$, $1\le i\le 5$,
and $\F\Gamma_{2,4}(14)$ of totally real algebraic number fields
were defined. The set $\F L^4$ consists of all ground fields of arithmetic
Lann\'er diagrams with $\ge 4$ vertices and consists of three fields of
degree $\le 2$. The set $\F T$ consists of all ground fields of arithmetic
triangles (plane) and has 13 fields of degree $\le 5$ (it includes $\F L^4$).
The set $\F\Gamma^{(4)}_i(14)$, $1\le i\le 5$, consists of all
ground fields of V-arithmetic finite edge polyhedra of minimality $14$ with
connected Gram graph having 4 vertices. They are determined by $5$ types
of graphs $\Gamma^{(4)}_i(14)$, $i=1,2,3,4,5$. The degrees of fields from
these sets are bounded by 22, 39, 53, 56, 54 respectively.
The set $\F\Gamma_{2,4}(14)$ consists of all ground fields of arithmetic
quadrangles (plane) of minimality $14$. Their degrees are bounded by 11.

The following result was obtained in \cite[Theorem 4.5]{Nik6}
using methods of \cite{Nik1} and \cite{Nik2} and results by Borel
\cite{Borel} and Takeuchi \cite{Tak4}.

\begin{theorem} (\cite{Nik6})
In dimensions $n\ge 6$, the ground field of
any arithmetic hyperbolic reflection group belongs to one of finite
sets of fields $\F L^4$, $\F T$, $\F\Gamma^{(4)}_i(14)$, $1\le i\le 5$,
and $\F\Gamma_{2,4}(14)$. In particular, its degree is bounded by $56$.
\label{thfor6}
\end{theorem}

In \cite[Secs 2 and 3]{Nik7}, further explicit finite sets of
fields $\F\Gamma^{(6)}_1(14)$, $\F\Gamma^{(6)}_2(14)$,
$\F\Gamma^{(6)}_3(14)$, $\F\Gamma^{(7)}_1(14)$,
$\F\Gamma^{(7)}_2(14)$, $\F \Gamma_{2,5}(14)$ were defined. The
sets $\F\Gamma^{(6)}_1(14)$, $\F\Gamma^{(6)}_2(14)$,
$\F\Gamma^{(6)}_3(14)$, $\F\Gamma^{(7)}_1(14)$,
$\F\Gamma^{(7)}_2(14)$ are defined by some V-arithmetic pentagon
graphs of minimality $14$. They are related to some fundamental
pentagons on hyperbolic plane. The degree of fields from
$\F\Gamma^{(6)}_1(14)$ is bounded by $56$, from
$\F\Gamma^{(6)}_2(14)$ by $75$, from $\F\Gamma^{(6)}_3(14)$ by
$138$, from $\F\Gamma^{(7)}_1(14)$ by $38$, from
$\F\Gamma^{(7)}_2(14)$ by $138$. The set $\F\Gamma_{2,5}(14)$
consists of all ground fields of arithmetic pentagons (plane) of
minimality $14$. Their degrees are bounded by 12.

Using similar, but much more complicated considerations,
we proved in \cite{Nik7}

\begin{theorem} In dimensions $n\ge 4$, the ground field of
any arithmetic hyperbolic reflection group belongs to one of finite
sets of fields $\F L^4$, $\F T$, $\F\Gamma^{(4)}_i(14)$, $1\le i\le 5$,
$\F\Gamma_{2,4}(14)$ and
$\F\Gamma^{(6)}_1(14)$,
$\F\Gamma^{(6)}_2(14)$,
$\F\Gamma^{(6)}_3(14)$, $\F\Gamma^{(7)}_1(14)$,\linebreak
$\F\Gamma^{(7)}_2(14)$,
$\F \Gamma_{2,5}(14)$.

In particular, its degree is bounded by $138$.
\label{thfor45}
\end{theorem}

Applying the same methods, here we want to extend this result to
$n\ge 3$, also considering $n=3$.

First, we introduce some other explicit finite sets of fields.
All of them are related to fundamental polygons on
hyperbolic plane.

We consider arithmetic reflection groups on hyperbolic plane with
fundamental polygons $\M_2 $ of minimality $14$. It means that
$\delta_1\cdot \delta_2<14$ for any $\delta_1,\delta_2\in P(\M_2)$.
The corresponding polygons $\M_2$ are also called arithmetic polygons of
the minimality $14$.

\begin{definition} We denote by $\Gamma_{2}(14)$ the set of
Gram graphs $\Gamma (P(\M_2))$ of all arithmetic polygons $\M_2$
of minimality $14$. The set $\F\Gamma_2(14)$ consists of all their
ground fields.
\label{F{2}(14)}
\end{definition}

It follows from results of Long, Maclachlan and Reid \cite{LMR}, Borel
\cite{Borel} and Takeuchi \cite{Tak4} (see \cite[Sec. 4.5]{Nik6})
that the degree of ground fields of arithmetic hyperbolic reflection
groups of dimension two is bounded by 44. Thus, the degree of fields
from $\F \Gamma_2(14)$ is also bounded by 44.

\begin{figure}
\begin{center}
\includegraphics[width=10cm]{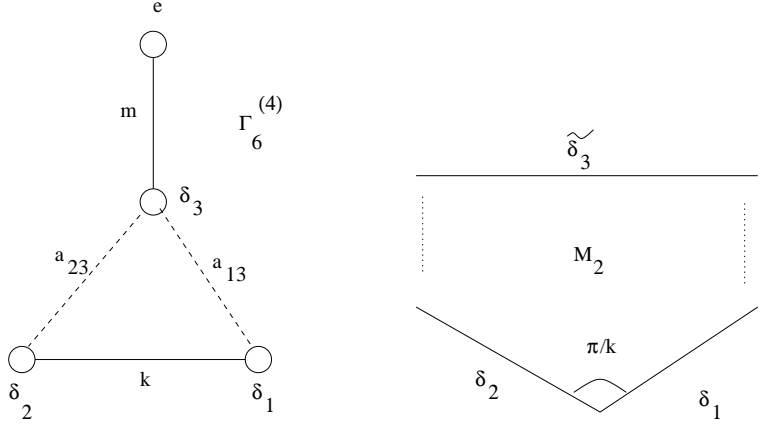}
\end{center}
\caption{3-dimensional graph $\Gamma_6^{(4)}$}
\label{3graph64}
\end{figure}

Let us consider V-arithmetic 3-dimensional chambers
which are defined by the Gram graphs $\Gamma^{(4)}_6$ with 4
vertices $\delta_1,\delta_2,\delta_3, e$
shown in Figure \ref{3graph64}.
It follows that the corresponding V-arithmetic chamber $\M$ satisfies
the following condition: the 2-dimensional face $\M_e$ of $\M$ which is
perpendicular to $e$ is an open fundamental triangle $\M_2$ bounded by
three lines perpendicular to
$$
P(\M_2)=\{\delta_1,\delta_2,\widetilde{\delta}_3\}
$$
for
$$
\widetilde{\delta}_3=\frac{\delta_3+\cos{(\pi/m)}e}{\sin{(\pi/m)}}\
$$
(see Figure \ref{3graph64}).
It has one angle $\pi/k$, $k\ge 3$,
defined by $\delta_1,\delta_2$. All its other sides don't intersect.
Planes $\Hh_{\delta_1}$ and $\Hh_{\delta_2}$ are perpendicular to $\Hh_e$,
and $\Hh_{\delta_3}$ has angle $\pi/m$, $m\ge 3$, with the plane $\Hh_e$.

\begin{definition} We denote by $\Gamma^{(4)}_6(14)$ the set of
all V-arithmetic 3-dimensional graphs $\Gamma^{(4)}_6$ (or the corresponding
3-dimensional V-arithmetic chambers) of minimality $14$.
Thus, inequalities $2<a_{ij}<14$ satisfy.
We denote by $\F \Gamma^{(4)}_6(14)$ the set of all their ground fields.
\label{def3graph64fields}
\end{definition}

\begin{figure}
\begin{center}
\includegraphics[width=12cm]{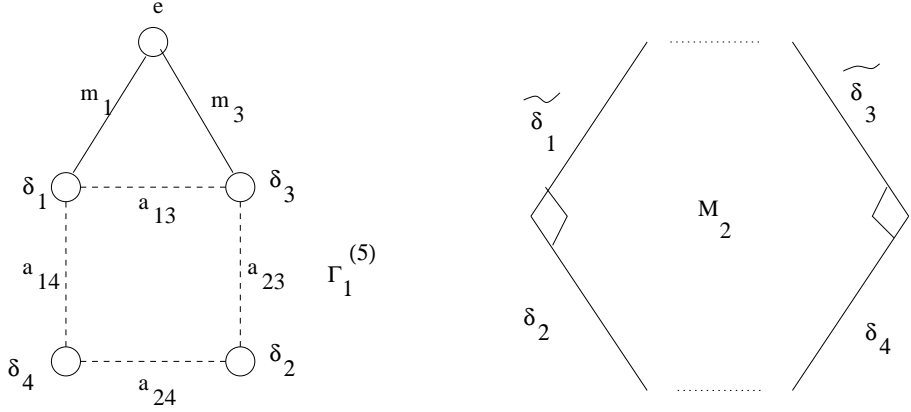}
\end{center}
\caption{3-dimensional graph $\Gamma_1^{(5)}$}
\label{3graph15}
\end{figure}

Let us consider V-arithmetic 3-dimensional chambers
which are defined by the Gram graphs $\Gamma^{(5)}_1$ with 5
vertices $\delta_1,\delta_2,\delta_3,\delta_4, e$
shown in Figure \ref{3graph15}.
It follows that the corresponding V-arithmetic chamber $\M$ satisfies
the following condition: the 2-dimensional face $\M_e$ of $\M$ which is
perpendicular to $e$ is an open fundamental quadrangle $\M_2$ bounded by
four lines perpendicular to
$$
P(\M_2)=\{\widetilde{\delta}_1,\delta_2,\widetilde{\delta}_3,\delta_4\}
$$
for
$$
\widetilde{\delta}_1=\frac{\delta_1+\cos{(\pi/m_1)}e}{\sin{(\pi/m_1)}},\ \
\widetilde{\delta}_3=\frac{\delta_3+\cos{(\pi/m_3)}e}{\sin{(\pi/m_3)}}
$$
(see Figure \ref{3graph15}).
It has two right angles defined by $\widetilde{\delta_1},\delta_2$ and
$\widetilde{\delta}_3,\delta_4$. All its other sides don't intersect.
Planes $\Hh_{\delta_2}$ and $\Hh_{\delta_4}$ are perpendicular to $\Hh_e$,
the planes $\Hh_{\delta_1}$ and $\Hh_{\delta_3}$ have angles $\pi/m_1$,
$m_1\ge 3$, and $\pi/m_3$, $m_3\ge 3$, with
the plane $\Hh_e$ respectively.

\begin{definition} We denote by $\Gamma^{(5)}_1(14)$ the set of
all V-arithmetic 3-dimensional graphs $\Gamma^{(5)}_1$ (or the corresponding
3-dimensional V-arithmetic chambers) of minimality $14$.
Thus, inequalities $2<a_{ij}<14$ satisfy.
We denote by $\F \Gamma^{(5)}_1(14)$ the set of all their ground fields.
\label{def3graph15fields}
\end{definition}

\begin{figure}
\begin{center}
\includegraphics[width=12cm]{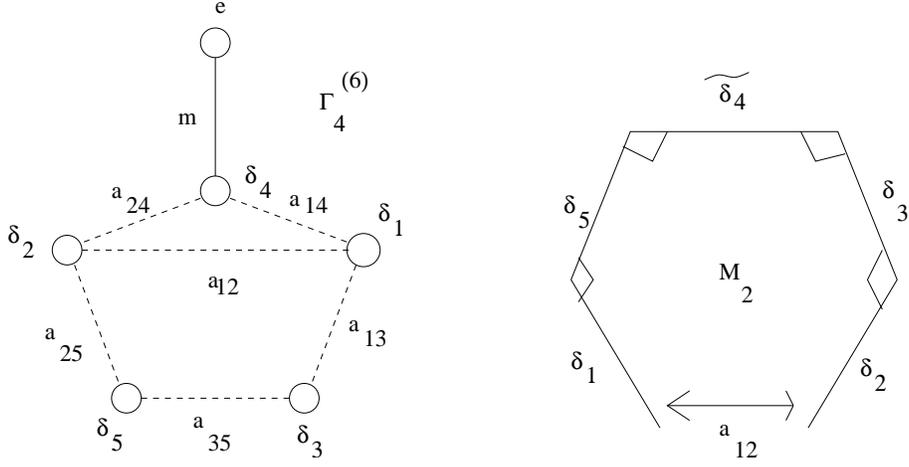}
\end{center}
\caption{3-dimensional graph $\Gamma_4^{(6)}$}
\label{3graph46}
\end{figure}

Let us consider V-arithmetic 3-dimensional chambers
which are defined by the Gram graphs $\Gamma^{(6)}_4$ with 6
vertices $\delta_1,\delta_2,\delta_3,\delta_4,\delta_5, e$
shown in Figure \ref{3graph46}.
It follows that the corresponding V-arithmetic chamber $\M$ satisfies
the following condition: the 2-dimensional face $\M_e$ of $\M$ which is
perpendicular to $e$ is an open fundamental pentagon $\M_2$ bounded by
5 lines perpendicular to
$$
P(\M_2)=\{\delta_1,\delta_2,\delta_3, \widetilde{\delta}_4,\delta_5\}
$$
for
$$
\widetilde{\delta}_4=\frac{\delta_1+\cos{(\pi/m)}e}{\sin{(\pi/m)}}
$$
(see Figure \ref{3graph46}).
It has four consecutive right angles defined by
$\delta_2$, $\delta_3$, $\widetilde{\delta_4}$, $\delta_5$, $\delta_1$
respectively. Its two consecutive sides perpendicular to $\delta_1$
and $\delta_2$ don't intersect.
All planes $\Hh_{\delta_i}$ are perpendicular to $\Hh_e$ except
the plane $\Hh_{\delta_4}$ which has the angle $\pi/m$ with
the plane $\Hh_e$.

\begin{definition} We denote by $\Gamma^{(6)}_4(14)$ the set of
all V-arithmetic 3-dimensional graphs $\Gamma^{(6)}_4$ (or the corresponding
3-dimensional V-arithmetic chambers) of minimality $14$.
Thus, inequalities $2<a_{ij}<14$ satisfy.
We denote by $\F \Gamma^{(6)}_4(14)$ the set of all their ground fields.
\label{def3graph46fields}
\end{definition}

We have

\begin{theorem} The sets of V-arithmetic graphs
$\Gamma_6^{(4)}(14)$, $\Gamma_1^{(5)}(14)$ and $\Gamma_4^{(6)}(14)$
are finite.

Degree of any field from $\F\Gamma_6^{(4)}(14)$ is bounded by $56$.

Degree of any field from $\F\Gamma_1^{(5)}(14)$ is bounded by
$909$.

Degree of any field from $\F \Gamma_4^{(6)}(14)$ is bounded by
$99$. \label{th3graphs}
\end{theorem}

\begin{proof} The proof requires long considerations and calculations.
It will be given in a special Section \ref{sec:3graphs}.
\end{proof}

We have the following main result of the paper.

\begin{theorem} In dimensions $n\ge 3$, the ground field of
any arithmetic hyperbolic reflection group belongs to one of finite
sets of fields $\F L^4$, $\F T$, $\F\Gamma^{(4)}_i(14)$, $1\le i\le 5$,
$\F\Gamma_{2,4}(14)$ (fields in dimension $\ge 6$), and
$\F\Gamma^{(6)}_1(14)$,
$\F\Gamma^{(6)}_2(14)$,
$\F\Gamma^{(6)}_3(14)$, $\F\Gamma^{(7)}_1(14)$,
$\F\Gamma^{(7)}_2(14)$,
$\F \Gamma_{2,5}(14)$ (additional fields for the dimensions $4$ and $5$),
and $\F\Gamma_6^{(4)}(14)$, $\F\Gamma_1^{(5)}(14)$,
$\F \Gamma_4^{(6)}(14)$, $\F \Gamma_2(14)$
(additional fields for the dimension $3$).

In particular, its degree is bounded by $909$. \label{th3fields}
\end{theorem}

\begin{proof} By \cite{Nik7}, if $n\ge 4$, the ground field $\bk$
belongs to one of sets
$\F L^4$, $\F T$, $\F\Gamma^{(4)}_i(14)$, $1\le i\le 5$,
$\F\Gamma_{2,4}(14)$, $\F\Gamma^{(6)}_1(14)$,
$\F\Gamma^{(6)}_2(14)$, $\F\Gamma^{(6)}_3(14)$, $\F\Gamma^{(7)}_1(14)$,
$\F\Gamma^{(7)}_2(14)$, $\F \Gamma_{2,5}(14)$.
Thus, further we can assume that the ground field $\bk$ does not belong
to these sets, and the dimension is equal to $n=3$.

Let $W$ be an arithmetic hyperbolic reflection group of dimension $n=3$,
$\M$ is its fundamental chamber, and $P(\M)$ is the set of all
vectors with square $-2$ which are perpendicular to codimension one faces
of $\M$ and directed outward. For $\delta \in P(\M)$ we denote by
$\Hh_\delta$ and  $\M_\delta$ the hyperplane and the codimension one face
$\M\cap \Hh_\delta$ respectively which is perpendicular to $\delta$.

By \cite{Nik1}, there exists $e\in P(\M)$ which defines
a narrow  face $\M_e$ of $\M$ of minimality $14$.
It means that $\delta_1\cdot \delta_2<14$ for any
$\delta_1,\delta_2\in P(\M,e)\subset P(\M)$. Here
$$
P(\M,e)=\{\delta \in P(\M)\ |\ \Hh_\delta\cap \Hh_e\not=\emptyset\}.
$$

Let us assume that $e\cdot \delta=0$ for any
$\delta\in P(\M,e)-\{e\}$, equivalently, all neighbouring 2-dimensional
faces of $\M$ to  the polygon $\M_e$ are perpendicular to $\M_e$.
Then $\M_e$ has $P(\M_e)=P(\M,e)-\{e\}$, and $\M_e$ is the fundamental
polygon for arithmetic hyperbolic plane reflection group with the same
ground field $\bk$ as for $W$. Since $P(\M_e)\subset P(\M)$,
it has minimality $14$. Then $\bk\in \F\Gamma_2(14)$ as required.

Thus, further we assume that $e\cdot \delta=2\cos(\pi/m)>0$ for one of
$\delta\in P(\M,e)-\{e\}$, equivalently, one of neighbouring to
the polygon $\M_e$ two-dimensional faces $\M_\delta$ of $\M$ is not
perpendicular to $\M_e$.

All cases when the polygon $\M_e$ has less than 6 sides were considered
in \cite{Nik6} and \cite{Nik7}. It was shown that then the ground field
$\bk$ belongs to one of sets of fields
$\F L^4$, $\F T$, $\F\Gamma^{(4)}_i(14)$, $1\le i\le 5$,
$\F\Gamma_{2,4}(14)$,
$\F\Gamma^{(6)}_1(14)$,
$\F\Gamma^{(6)}_2(14)$,
$\F\Gamma^{(6)}_3(14)$, $\F\Gamma^{(7)}_1(14)$,
$\F\Gamma^{(7)}_2(14)$,
$\F \Gamma_{2,5}(14)$. Thus, further we additionally
assume that $\M_e$ has at least 6 sides and the ground field $\bk$
does not belong to any of these sets of fields.

By \cite[Lemma 4.3]{Nik6}, if $\bk$ does not belong to $\F L^4$,
$\F T$ and $\F \Gamma_i^{(4)}(14)$, $1\le i\le 4$, then the
Coxeter graph $C(v)$ of any vertex $v\in \M_e$ has all connected
components having only one or two vertices. If additionally $\bk$
does not belong to $\F \Gamma_5^{(4)}(14)$, then the hyperbolic
connected component of the edge graph $\Gamma(r)$ defined by any
edge $r=v_1v_2\subset \M_e$ has $\le 3$ vertices. Further we mark
these facts as (*).

By (*), both angles of $\M_e$ at the edge $\M_{e,\delta}$
perpendicular to $\delta$ are right. Moreover, if $f_1,f_2\in P(\M,e)$
define two neighbouring edges of the edge $\M_{e,\delta}$ of $\M_e$,
then $e\cdot f_1=e\cdot f_2=\delta\cdot f_1=\delta\cdot f_2=0$.

Assume that the polygon $\M_e$ has a non-right angle with edges
perpendicular to $\delta_1,\delta_2\in P(\M, e)$. By (*),  then
$\delta_1\cdot\delta_2=2\cos(\pi/k)>0$, and $\{\delta_1,\delta_2\}$ are
perpendicular to $\{e,\delta\}$. Then $\delta_1,\delta_2, \delta_3=\delta$
and $e$ have Gram graph $\Gamma_6^{(4)}(14)$, and the ground $\bk$ belongs to
$\F \Gamma_6^{(4)}(14)$ as required.

Now assume that all angles of the polygon $\M_e$ are right and there exist two
elements $\delta_1,\delta_3\in P(\M,e)-\{e\}$ such that
$e\cdot \delta_1=2\cos(\pi/m_1)$ and
$e\cdot \delta_3=2\cos(\pi/m_3)$. By (*), then $\delta_1\cdot \delta_3>0$,
and  $\delta_1,\delta_3$ are perpendicular to two not consecutive edges
of the polygon $\M_e$.
Since $\M_e$ has more than 5 vertices, we can find two their
neighbouring edges perpendicular to $\delta_2, \delta_4\in P(\M,e)$
such that the Gram graph of $\delta_1,\delta_2, \delta_3, \delta_4,e$
is $\Gamma_1^{(5)}$. Then the ground field belongs to $\F\Gamma_1^{(5)}(14)$
as required.

Now assume that all angles of the polygon $\M_e$ are right and there exists
only one element $\delta_4\in P(\M,e)-\{e\}$ such that
$\delta_4\cdot e\not=0$. Then $\delta_4\cdot e=2\cos(\pi/m)>0$. Since $\M_e$
has at least 6 sides, we can find 5 consecutive sides of $\M_e$ perpendicular
to $\delta_1,\delta_5,\delta_4,\delta_3,\delta_2\in P(\M,e)$
such that Gram graph
of $\delta_1,\delta_2, \delta_3, \delta_4, \delta_5, e$ is $\Gamma_4^{(6)}$.
Then the ground field $\bk$ belongs to $\F \Gamma_4^{(6)}(14)$ as
required.

This finishes the proof of the theorem. \end{proof}



\section{V-arithmetic 3-graphs $\Gamma_6^{(4)}(14)$,
$\Gamma_1^{(5)}(14)$, $\Gamma_4^{(6)}(14)$  and their
fields}\label{sec:3graphs}

Here we prove Theorem \ref{th3graphs}.


\subsection{Some general results.}
\label{subsec:genresults} We use the following general results
from \cite{Nik2}.

\begin{theorem} (\cite[Theorem 1.2.1]{Nik2}) Let $\bff$ be a totally
real algebraic number field, and let each embedding
$\sigma:\bff\to \br$ corresponds to an interval
$[a_\sigma,b_\sigma]$ in $\br$ where
$$
\prod_{\sigma }{\frac{b_\sigma-a_\sigma}{4}}<1.
$$
In addition, let the natural number $m$ and the intervals
$[s_1,t_1],\dots, [s_m,t_m]$ in $\br$ be fixed. Then there exists
a constant $N(s_i,t_i)$ such that, if $\alpha$ is a totally real
algebraic integer and if the following inequalities hold for the
embeddings $\tau:\bff(\alpha) \to \br$:
$$
s_i\le \tau(\alpha)\le t_i\ \ for\ \ \tau=\tau_1,\dots ,\tau_m,
$$
$$
a_{\tau | \bff}\le \tau(\alpha)\le b_{\tau | \bff}\ \ for\ \
\tau\not=\tau_1,\dots,\tau_m,
$$
then
$$
[\bff(\alpha):\bff]\le N(s_i,t_i).
$$
\label{th121}
\end{theorem}

\begin{theorem} (\cite[Theorem 1.2.2]{Nik2})
Under the conditions of Theorem \ref{th121}, $N(s_i,t_i)$ can be
taken to be $N(s_i,t_i)=N_0$,
 where $N_0$ is the least natural number
solution of the inequality
\begin{equation}
N_0M\ln{(1/R)} - M\ln{(N_0+1)}-\ln{B}\ge \ln{S}. \label{cond for
n}
\end{equation}

Here
\begin{equation}
M=[\bff : \bq],\ \ \ B=2\sqrt{|{\rm discr\ } \bff|}; \label{MB}
\end{equation}
\begin{equation}
R=\sqrt{\prod_\sigma {\frac{b_\sigma-a_\sigma}{4}}},\ \ \
S=\prod_{i=1}^{m}{\frac{2er_i}{b_{\sigma_i}-a_{\sigma_i}}}
\label{RS}
\end{equation}
where
\begin{equation}
\sigma_i=\tau_i|\bff,\ \ \  r_i=\max\{{|b_i-a_{\sigma_i}|,
|b_{\sigma_i}-a_i|}\}. \label{ri}
\end{equation}
\label{th122}
\end{theorem}

We note that the proof of Theorems \ref{th121} and \ref{th122}
uses a variant of Fekete's Theorem (1923) about existence of
non-zero integer polynomials of bounded degree which differ only
slightly from zero on appropriate intervals. See \cite[Theorem
1.1.1]{Nik2}.

Below we will apply these results in two cases which are very
similar to used in \cite[Sec. 5.5]{Nik6} and \cite{Nik2}. Cases 1
and 2 below are natural generalizations of Cases 1 and 2 which
were considered in \cite[Sec. 3.1]{Nik7}.

\medskip

{\bf Case 1.} For a natural $l\ge 3$ we denote
$\bff_l=\bq(\cos{(2\pi/l)})$. We consider a totally real algebraic
number field $\bk$ where $\bff_l\subset \bk=\bq(\alpha)$, and the
algebraic integer $\alpha$ satisfies
\begin{equation}
-a_1\sigma(\sin^2{(\pi/l)})<\sigma(\alpha)<a_2\sigma(\sin^2{(\pi/l)})
\label{case1.cond1sigma}
\end{equation}
for all $\sigma:\bk \to \br$ such that $\sigma\not=\sigma^{(+)}$, and
\begin{equation}
b_1<\sigma^{(+)}(\alpha)<b_2
\label{case1.cond1sigmapl}
\end{equation}
where $\sigma^{(+)}:\bk\to \br$ is the identity.  We assume that
$a_1\ge 0$, $a_2\ge 0$ and $0<a=\max\{a_1,a_2\}<4$.
We assume that $b_1<b_2$ and
denote $b=\max\{|b_1|, |b_2|\}$. Also we assume that $a\le b$.
We want to estimate $[\bk:\bff_l]=N_0$ and
$N=[\bk:\bq]=N_0\cdot [\bff_l:\bq]$ from above.

For $l\ge 3$, we have $[\bff_l:\bq]=\varphi(l)/2$ where $\varphi(l)$ is
the Euler function, and
$N_{\bff_l/\bq}(\sin^2{(\pi/l)})=\gamma(l)/4^{\varphi(l)/2}$ where
\begin{equation}
N_{\bff_l/\bq }(4\sin^2{(\pi/l)})=\gamma(l)=\left\{
\begin{array}{cl}
p &\ {\rm if}\  l=p^t>2\ {\rm where }\ p\ {\rm is\ prime,} \\
1 &\ {\rm otherwise.}
\end{array}\
\label{normsin1}
 \right .
\end{equation}
We have
$$
\frac{ba^N|N_{\bk/\bq}(\sin^2{(\pi/l)})|}{a\sin^2{(\pi/l)}}>
|N_{\bk/\bq}(\alpha)|\ge 1
$$
and
$$
\frac{b(a/4)^N\gamma(l)^{2N/\varphi(l)}}{a\sin^2{(\pi/l)}}>1.
$$
Equivalently, we have
\begin{equation}
N\left(\ln{\frac{2}{\sqrt{a}}-\frac{\ln{\gamma(l)}}{\varphi(l)}}\right)
<\ln{\sqrt{\frac{b}{a}}}-\ln{\sin{\frac{\pi}{l}}}\ \ \text{and}\ \
(\varphi(l)/2)|N.
\label{case1.1}
\end{equation}
Since $\gamma (l)\le l$, $\varphi(l)\ge Cl/\ln(\ln{l})$ for $l\ge
6$ where $C=\varphi(6)\ln{(\ln{6})}/6\ge 0.194399$,
$\sin{(\pi/l)}\le \pi/l$ for $l\ge 3$, there exists only finite
number of $l\ge 3$ such that \eqref{case1.1} has solutions $N\in \bn$.

More exactly, there exists only finite number of {\it exceptional
$l\ge 3$} such that
\begin{equation}
\ln{\frac{2}{\sqrt{a}}-\frac{\ln{\gamma(l)}}{\varphi(l)}}\le 0.
\label{case1.2}
\end{equation}
All non-exceptional $l$ satisfy the inequality
\begin{equation}
(\varphi(l)/2)\left(\ln{\frac{2}{\sqrt{a}}-
\frac{\ln{\gamma(l)}}{\varphi(l)}}\right)
<\ln{\sqrt{\frac{b}{a}}}-\ln{\sin{\frac{\pi}{l}}}.
\label{case1.3}
\end{equation}
Remark that exceptional $l$ also satisfy this inequality.

If $\gamma(l)=1$, \eqref{case1.3} implies that $l$ satisfies the
inequality
\begin{equation}
 {(C/2)\ln{(2/\sqrt{a})}}\,l<{\left(\ln{l}+\ln{(\sqrt{(b/a)}/\pi)}\right)
\ln{\ln{l}}}\,.
\label{case1.4}
\end{equation}
It follows that
\begin{equation}
l<L_0 \label{case1.5}
\end{equation}
where $L_0>3$ satisfies
\begin{equation}
{(C/2)\ln{(2/\sqrt{a})}}\,L_0\ge
{\left(\ln{L_0}+\ln{(\sqrt{(b/a)}/\pi)}\right) \ln{\ln{L_0}}}\,.
\label{case1.6}
\end{equation}
If $l=p^t$ where $p$ is prime, \eqref{case1.3} implies that $l$
satisfies the inequality
\begin{equation}
(C/2)\Delta(a)\,l<{\left(\ln{l}+\ln{(\sqrt{(b/a)}/\pi)}\right)
\ln{\ln{l}}}\,
\label{case1.7}
\end{equation}
where
\begin{equation}
\Delta(a)=\min_{l=p^t\ge L_0}
\left\{\ln{\frac{2}{\sqrt{a}}-\frac{\ln{\gamma(l)}}{\varphi(l)}>0}\right\}.
\label{case1.8}
\end{equation}
It follows that
\begin{equation} l<L_1
\end{equation}
where $L_1\ge L_0$ is a solution of the inequality
\begin{equation}
(C/2)\Delta(a)\,L_1\ge
{\left(\ln{L_1}+\ln{(\sqrt{(b/a)}/\pi)}\right)
\ln{\ln{L_1}}}\,.
\label{case1.9}
\end{equation}
Thus, to find all non-exceptional $l$ satisfying \eqref{case1.3},
we should check \eqref{case1.3} for all $l$ such that $3\le l<
L_1$, moreover, if $L_0\le l<L_1$, we can assume that $l=p^t$.
Their number is finite, and all of them can be effectively found.

For non-exceptional $l$ satisfying \eqref{case1.3}, we obtain
bounds
\begin{equation}
N_0=[\bk:\bff_l]\le \left[
\frac{\ln{\sqrt{b/a}}-\ln{\sin{(\pi/l)}}}
{(\varphi(l)/2)\left(\ln{(2/\sqrt{a})}-(\ln{\gamma(l)})/\varphi(l)\right)}
\right]
\label{case1.10}
\end{equation}
and
\begin{equation}
N=[\bk:\bq]\le \left[ \frac{\ln{\sqrt{b/a}}-\ln{\sin{(\pi/l)}}}
{(\varphi(l)/2)\left(\ln{(2/\sqrt{a})}-(\ln{\gamma(l)})/\varphi(l)\right)}
\right]\cdot(\varphi(l)/2)\,.
\label{case1.11}
\end{equation}
This using of the norm, we call the {\it Method B} (like in
\cite[Sec. 5.5]{Nik6}).

On the other hand, for fixed $l$, we obtain a bound for $N_0$
using Theorems \ref{th121} and \ref{th122} applied to 
$\bff=\bff_l$ and $\alpha$. We can take
\begin{equation}
R=\sqrt{|N_{\bff_l/\bq}(\sin^2{(\pi/l)})|
\left((a_1+a_2)/4\right)^{\varphi(l)/2}}=
\left(\frac{\gamma(l)^{1/\varphi(l)}(a_1+a_2)^{1/2}}{4}\right)^{\varphi(l)/2},
\label{case1.12}
\end{equation}
where
\begin{equation}
R<1\ \text{ if and only if }\
\ln{\frac{4}{\sqrt{a_1+a_2}}}-\frac{\ln{\gamma(l)}}{\varphi(l)}>0\,,
\label{case1.12a}
\end{equation}
\begin{equation}
M=[\bff_l:\bq]=\frac{\varphi(l)}{2},\ \
B=2\sqrt{|\text{discr}\,\bff_l|}
\label{case1.13}
\end{equation}
where the discriminant $|\text{discr}\,\bff_l|$ is given in
\eqref{discrFl}, and
\begin{equation}
S=\frac{2e\max\{a_2,b_2,a_2-b_1,a_1,-b_1,b_2+a_1\}}{(a_1+a_2)\sin^2{(\pi/l)}}.
\label{case1.14}
\end{equation}
Then $[\bk:\bff_l]\le n_0$ and $[\bk:\bq]\le n_0\varphi(l)/2$
where $n_0$ is the least natural solution of the inequality
\eqref{cond for n}
\begin{equation} n_0M\ln{(1/R)} -
M\ln{(n_0+1)}-\ln{B}\ge \ln{S}.
\label{case1.15}
\end{equation}
In particular, this gives a bound for $[\bk:\bq]$ for exceptional
$l$ satisfying \eqref{case1.12a} and improves the bound
\eqref{case1.10} for $N_0$ when it is poor, which also improves
the bound for $[\bk:\bq]$. This using of Theorems \ref{th121},
\ref{th122}, we call the {\it Method A} (like in \cite[Sec.
5.5]{Nik6}).

We shall apply these Methods A and B to $\Gamma_4^{(6)}(14)$ in Sect.
\ref{subsec:3graph46}.

\medskip


{\bf Case 2.} For natural $k\ge s\ge 3$, we denote
$\bff_{k,s}=\bq(\cos{(2\pi/k)},\,\cos{(2\pi/s)})$. We consider a
totally real algebraic number field $\bk$ where $\bff_{k,s}\subset
\bk=\bq(\alpha)$, and the algebraic integer $\alpha$ satisfies
\begin{equation}
-a_1\sigma(\sin^2{(\pi/k)}\sin^2{(\pi/s)})<\sigma(\alpha)
<a_2\sigma(\sin^2{(\pi/k)}\sin^2{(\pi/s)})
\label{case2.cond1sigma}
\end{equation}
for all $\sigma:\bk \to \br$ such that $\sigma\not=\sigma^{(+)}$, and
\begin{equation}
b_1<\sigma^{(+)}(\alpha)<b_2
\label{case2.cond1sigmapl}
\end{equation}
where $\sigma^{(+)}:\bk\to \br$ is the identity.  We assume that
$a_1\ge 0$, $a_2\ge 0$ and
$0<a=\max\{a_1,a_2\}<16$. We assume that  $b_1<b_2$ and denote
$b=\max\{|b_1|, |b_2|\}$. Also we assume that $a\le b$.
We want to estimate $[\bk:\bff_{k,s}]=N_0$
and $N=[\bk:\bq]=N_0[\bff_{k,s}:\bq]$ for non-exceptional $k$ and
$s$ where {\it $l\ge 3$ is called exceptional} if
\begin{equation}
\ln{\frac{4}{\sqrt{a}}}-\frac{\ln{\gamma(l)}}{\varphi(l)}\le 0.
\label{case2.excepl}
\end{equation}
We also assume that $k\ge s\ge s_0\ge 3$ where $s_0\ge 3$ is fixed.

We have $[\bff_{k,s}:\bq]=\varphi([k,s])/2\rho(k,s)$ where
$\rho(k,s)=1$ or $2$ is given in \eqref{defrho}, and
$N_{\bff_l/\bq}(\sin^2{(\pi/l)})=\gamma(l)/4^{\varphi(l)/2}$ where
$\gamma(l)$ is given in \eqref{normsin1}. We have
$$
\frac{ba^N|N_{\bk/\bq}(\sin^2{(\pi/k)})|
|N_{\bk/\bq}(\sin^2{(\pi/s)})| }
{a\sin^2{(\pi/k)}\sin^2{(\pi/s)}}> |N_{\bk/\bq}(\alpha)|\ge 1
$$
and
$$
\frac{b(a/16)^N\gamma(k)^{2N/\varphi(k)}\gamma(s)^{2N/\varphi(s)}}
{a\sin^2{(\pi/k)}\sin^2{(\pi/s)}}>1.
$$
Equivalently, we obtain
\begin{equation}
N\left(\ln{\frac{4}{\sqrt{a}}-\frac{\ln{\gamma(k)}}{\varphi(k)}-
\frac{\ln{\gamma(s)}}{\varphi(s)} }\right)
<\ln{\sqrt{\frac{b}{a}}}-\ln{\sin{\frac{\pi}{k}}}-
\ln{\sin{\frac{\pi}{s}}}\ \text{and}\
\frac{\varphi([k,s])}{2\rho(k,s)}\left|\right. N.
\label{case2.1}
\end{equation}
Since $\gamma (l)\le l$, $\varphi(l)\ge Cl/\ln(\ln{l})$ for $l\ge
6$ where $C=\varphi(6)\ln{(\ln{6})}/6$, \linebreak
$\sin{(\pi/l)}\le \pi/l$
for $l\ge 3$, there exists only finite number of pairs $(k,s)$ such
that \eqref{case2.1} has solutions $N\in \bn$
where both $k$ and $s$ are non-exceptional.

More exactly, there exists only finite number of {\it exceptional
pairs $(k,s)$} where a pair $(k,s)$ (consisting of non-exceptional
$k$ and $s$) is called exceptional  if
\begin{equation}
\ln{\frac{4}{\sqrt{a}}}-\frac{\ln{\gamma(k)}}{\varphi(k)}-
\frac{\ln{\gamma(s)}}{\varphi(s)}\le 0.
\label{case2.2}
\end{equation}
All non-exceptional pairs $(k,s)$ satisfying \eqref{case2.1}
satisfy the inequality
\begin{equation}
\frac{\varphi([k,s])}{2\rho(k,s)}\cdot \left(
\ln{\frac{4}{\sqrt{a}}}-\frac{\ln{\gamma(k)}}{\varphi(k)}-
\frac{\ln{\gamma(s)}}{\varphi(s)}\right) <
\ln{\sqrt{\frac{b}{a}}}-\ln{\sin{\frac{\pi}{k}}}-
\ln{\sin{\frac{\pi}{s}}}\ .
\label{case2.3}
\end{equation}
Remark that exceptional pairs $(k,s)$ also satisfy this inequality.

If $\gamma(k)=\gamma(s)=1$ and $k\ge s$, \eqref{case2.3} implies
\begin{equation}
{(C/2)\ln{(4/\sqrt{a})}}
k<{\left(2\ln{k}+\ln{(\sqrt{(b/a)}/\pi^2)}\right)\ln{\ln{k}}}.
\label{case2.4}
\end{equation}
It follows that
\begin{equation}
s_0\le s\le k< K_0
\label{case2.5}
\end{equation}
where $K_0>3$ satisfies
\begin{equation}
{(C/2)\ln{(4/\sqrt{a})}}
K_0\ge {\left(2\ln{K_0}+\ln{(\sqrt{(b/a)}/\pi^2)}\right)\ln{\ln{K_0}}}.
\label{case2.6}
\end{equation}
If one of $\gamma(k)$, $\gamma(s)$ is not equal to $1$, then
\eqref{case2.3} implies for non-exceptional pairs $(k,s)$ that
\begin{equation}
{(C/2)\Delta_1(a)}
k<{\left(2\ln{k}+\ln{(\sqrt{(b/a)}/\pi^2)}\right)\ln{\ln{k}}}
\label{case2.7}
\end{equation}
where
\begin{equation}
\Delta_1(a)=\min_{k\ge s\ge s_0, k\ge K_0}\left\{\ln{\frac{4}{\sqrt{a}}}-
\frac{\ln{\gamma(s)}}{\varphi(s)}-
\frac{\ln{\gamma(k)}}{\varphi(k)}\,>\,0\right\}\,.
\label{case2.8}
\end{equation}
It follows that
\begin{equation}
s_0\le s\le k<K_1
\label{case2.9}
\end{equation}
where $K_1\ge K_0$ is a solution of the inequality
\begin{equation}
{(C/2)\Delta_1(a)}
K_1\ge {\left(2\ln{K_1}+\ln{(\sqrt{(b/a)}/\pi^2)}\right)\ln{\ln{K_1}}}.
\label{case2.10}
\end{equation}
Thus, to find all non-exceptional pairs $(k,s)$ satisfying \eqref{case2.3},
we should check \eqref{case2.3} for all $s_0\le s\le k<K_1$; moreover, if
$K_0\le k\le K_1$, we can assume that one of $k$ and $s$ is equal to $p^t$
where $p$ is prime. The number of such pairs is finite, and all of them
can be effectively found.

For such non-exceptional pairs $(k,s)$ satisfying \eqref{case2.3},
we obtain bounds
\begin{equation}
N_0=[\bk:\bff_{k,s}]\le \left[\frac
{\ln{\sqrt{\frac{b}{a}}}-\ln{\sin{\frac{\pi}{k}}}-
\ln{\sin{\frac{\pi}{s}}}} {\frac{\varphi([k,s])}{2\rho(k,s)}\cdot
\left( \ln{\frac{4}{\sqrt{a}}}-\frac{\ln{\gamma(k)}}{\varphi(k)}-
\frac{\ln{\gamma(s)}}{\varphi(s)}\right)}\right]
\label{case2.11}
\end{equation}
and
\begin{equation}
N=[\bk:\bq]\le \left[\frac
{\ln{\sqrt{\frac{b}{a}}}-\ln{\sin{\frac{\pi}{k}}}-
\ln{\sin{\frac{\pi}{s}}}} {\frac{\varphi([k,s])}{2\rho(k,s)}\cdot
\left( \ln{\frac{4}{\sqrt{a}}}-\frac{\ln{\gamma(k)}}{\varphi(k)}-
\frac{\ln{\gamma(s)}}{\varphi(s)}\right)}\right]\cdot
\frac{\varphi([k,s])}{2\rho(k,s)}\ .
\label{case2.12}
\end{equation}
This using of the norm, we call the {\it Method B} (like in
\cite[Sec. 5.5]{Nik6}).

On the other hand, for a fixed pair $(k,s)$, we can obtain a bound
for $N_0$ using Theorems \ref{th121} and \ref{th122} applied 
to $\bff=\bff_{k,s}$ and $\alpha$. We can take
$$
R=\sqrt{|N_{\bff_{k,s}/\bq}(\sin^2{\frac{\pi}{k}}
\sin^2{\frac{\pi}{s}})|
\left(\frac{a_1+a_2}{4}\right)^{\frac{\varphi([k,s])}{2\rho(k,s)}}}=
$$
\begin{equation}
=\left(\frac{\gamma(k)^{\frac{1}{\varphi(k)}}
\gamma(s)^{\frac{1}{\varphi(s)}}
(a_1+a_2)^{\frac{1}{2}}}{8}\right)^{\frac{\varphi([k,s])}{2\rho(k,s)}}
\label{case2.13}
\end{equation}
where
\begin{equation}
R<1\ \text{if and only if}\
\ln{\frac{8}{\sqrt{a_1+a_2}}}-\frac{\ln{\gamma(k)}}{\varphi(k)}-
\frac{\ln{\gamma(s)}}{\varphi(s)}>0,
\label{case2.14}
\end{equation}
\begin{equation}
M=[\bff_{k,s}:\bq]=\frac{\varphi([k,s])}{2\rho(k,s)},\ \
B=2\sqrt{|\text{discr}\,\bff_{k,s}|}
\label{case2.15}
\end{equation}
where the discriminant $|\text{discr}\,\bff_{k,s}|$ is given in
\eqref{discrFks1} and \eqref{discrFks2}, and
\begin{equation}
S=\frac{2e\max\{a_2,b_2,a_2-b_1,a_1,-b_1,b_2+a_1\}}
{(a_1+a_2)\sin^2{(\pi/s)}\sin^2{(\pi/k)}}.
\label{case2.16}
\end{equation}
For all pairs $(k,s)$ satisfying \eqref{case2.14}, we obtain the
bounds $[\bk:\bff_{k,s}]\le n_0$ and $[\bk:\bq]\le
n_0\varphi([k,s])/(2\rho(k,s))$ where $n_0$ is the least natural
solution of the inequality \eqref{cond for n},
$$
n_0M\ln{(1/R)} - M\ln{(n_0+1)}-\ln{B}\ge \ln{S}.
$$
For $a<16$ and $k,s\ge 3$, all
pairs $(k,s)$, except finite number, satisfy \eqref{case2.14},
and we can apply this method to all these pairs.
In particular, this gives a bound for $[\bk:\bq]$ for all exceptional pairs
$(k,s)$ satisfying \eqref{case2.14}, and it improves the bound
\eqref{case2.11} for $N_0$ when it is poor,
which also improves the bound for $[\bk:\bq]$.
This using of Theorems \ref{th121}, \ref{th122}, we call the {\it
Method A} (like in \cite[Sec. 5.5]{Nik6}).

We apply these Methods A and B to $\Gamma_6^{(4)}(14)$ in Sect.
\ref{subsec:3graph64} and to
$\Gamma_1^{(5)}(14)$ in Sect. \ref{subsec:3graph15}.


\subsection{V-arithmetic 3-graphs $\Gamma_6^{(4)}(14)$
and their fields.}\label{subsec:3graph64}

Here we consider V-arithmetic 3-dimensional graphs
$\Gamma_6^{(4)}(14)$ (see Figure \ref{3graph64}) and their fields.

This case had been considered in \cite[Sec. 3.1]{Nik7} where the
bound for degrees of fields $\bk$ from $\F \Gamma_1^{(6)}(14)$ was
obtained. To get this bound, in \cite[Sec. 3.2]{Nik7} the bound
for fields defined by the subgraph $\Gamma_6^{(4)}(14)$ of this
graph was obtained. This uses methods A and B of Case 2 in Sect.
\ref{subsec:genresults} applied to $a_1=0$, $a_2=4$, $b_1=12$ and
$b_2=28^2$. The bound is $[\bk:\bq]\le 56$.

\subsection{V-arithmetic 3-graphs $\Gamma_1^{(5)}(14)$
and their fields.}\label{subsec:3graph15}

Here we consider V-arithmetic 3-dimensional graphs
$\Gamma_1^{(5)}(14)$ (see Figure \ref{3graph15}) and their fields.

\begin{figure}
\begin{center}
\includegraphics[width=11cm]{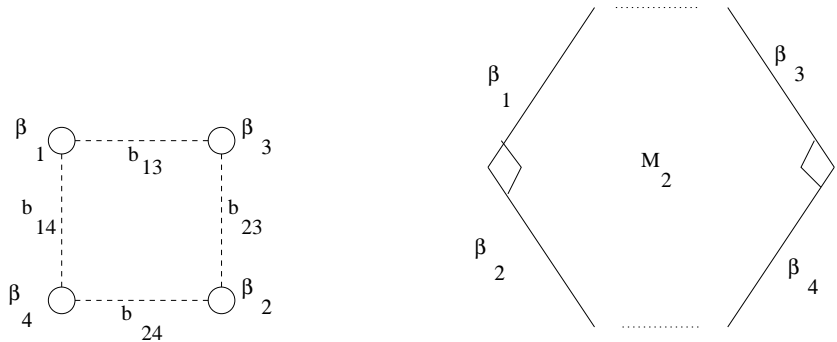}
\end{center}
\caption{The graph of $\M_2$ for $\Gamma_1^{(5)}$}
\label{plgraph15}
\end{figure}

First, let us consider the corresponding plane graph defined by
$\beta_1=\widetilde{\delta_1}$, $\beta_2=\delta_2$,
$\beta_3=\widetilde{\delta}_3$, $\beta_4=\delta_4$ which
give $P(\M_2)$.
We denote $b_{ij}=\beta_i\cdot \beta_j$ when it is not $0$.
This graph is given in Figure \ref{plgraph15}.

Any three elements from $\beta_1,\dots,\beta_4$ generate the form
defining the hyperbolic plane. Thus the determinant of their Gram
matrix must be positive for geometric embedding $\sigma^{(+)}$ and
must be negative for $\sigma\not=\sigma^{(+)}$. For example,
for $\beta_1, \beta_2, \beta_3$ it is equal to
$-8+2b_{13}^2+2b_{23}^2$. Thus, for $\sigma$ we obtain inequalities
$b_{13}^2+b_{23}^2<4$. Moreover, the determinant
$$
16+b_{13}^2 b_{24}^2+b_{14}^2
b_{23}^2-4b_{13}^2-4b_{14}^2-4b_{23}^2-4b_{24}^2-2b_{13}b_{14}b_{23}b_{24}
$$
of the Gram matrix of all four elements $\beta_1,\dots,\beta_4$ is
$0$. Combining all these conditions, we obtain the following
conditions on $\M_2$ for $\sigma\not=\sigma^{(+)}$:
\begin{equation}
\left\{
\begin{array}{l}
b_{13}b_{14}b_{23}b_{24}=8-2b_{13}^2-2b_{14}^2-2b_{23}^2-2b_{24}^2+
(b_{13}^2b_{24}^2+b_{14}^2b_{23}^2)/2\\
b_{13}^2+b_{23}^2<4\\
b_{23}^2+b_{24}^2<4\\
b_{24}^2+b_{14}^2<4\\
b_{14}^2+b_{13}^2<4\ .
\end{array}
\right.
\label{equationgam15}
\end{equation}
It is easy to find minimum and maximum of
$b_{13}b_{14}b_{23}b_{24}$ under the closure of these conditions
which shows that
\begin{equation}
-4< \sigma(b_{13}b_{14}b_{23}b_{24})\le 1. \label{ingam15}
\end{equation}
Here minimum is achieved for $b_{ij}=\pm \sqrt{2}$ where the
number of $(-)$ is odd, and maximum is achieved for $b_{ij}=\pm 1$
where the number of $(-)$ is even.
From expressions of $\beta_i$
using $\delta_i$ and $e$, we get
\begin{equation}
b_{13}b_{14}b_{23}b_{24}=\frac
{a_{13}a_{14}a_{23}a_{24}+
2\cos{\frac{\pi}{m_1}}\cos{\frac{\pi}{m_3}}a_{14}a_{23}a_{24}}
{\sin^2{\frac{\pi}{m_1}}{\sin^2{\frac{\pi}{m_3}}}}.
\label{bagam15}
\end{equation}
We consider the algebraic integer $\alpha\in \bk$ which is
\begin{equation}
\alpha= 2a_{13}a_{14}a_{23}a_{24}+
4\cos{\frac{\pi}{m_1}}\cos{\frac{\pi}{m_3}}a_{14}a_{23}a_{24}.
\label{alpha15}
\end{equation}
From \eqref{ingam15} and \eqref{bagam15}, we get
\begin{equation}
-8\sigma(\sin^2{\frac{\pi}{m_1}}\sin^2{\frac{\pi}{m_3}})<
\sigma(\alpha)\le
2\sigma(\sin^2{\frac{\pi}{m_1}}\sin^2{\frac{\pi}{m_3}}).
\label{alphasigma15}
\end{equation}
For the geometric embedding $\sigma^{(+)}$, we have
\begin{equation}
2\cdot 2^4+2^3=5\cdot 2^3<\sigma^{(+)}(\alpha)<2\cdot 14^4+4\cdot
14^3=32\cdot 14^3. \label{alphasigma+}
\end{equation}

It follows that $\bk=\bq(\alpha)$. Since
$8 <16$, this case is
similar to considered in \cite[Sec. 5.5]{Nik6}.


In this case, considering $\alpha$ from \eqref{alpha15}, we apply
the methods A and B of Case 2 in Sec. \ref{subsec:genresults} to
$a_1=8$, $a_2=2$ (then $a=8$), $b_1=5\cdot 2^3$, $b_2=32\cdot
14^3$ (then $b=32\cdot 14^3$) and $k:=\max\{m_1,m_3\}$,
$s:=\min\{m_1,m_3\}$ where $k\ge s\ge 3$, and $s_0=3$.

At first, we apply the Method B. Exceptional $l\ge 3$ satisfy
\eqref{case2.excepl} which is
\begin{equation}
\ln{\sqrt{2}}-\frac{\ln{\gamma(l)}}{\varphi(l)}\le 0.
\label{Gamma1.5.0}
\end{equation}
It follows that $l=3,4,5$ are the only exceptional.

All exceptional pairs $(k,s)$ where $k\ge s \ge 6$, that is when
\eqref{case2.2} which is
\begin{equation}
\ln{\sqrt{2}}-\frac{\ln{\gamma(k)}}{\varphi(k)}-
\frac{\ln{\gamma(s)}}{\varphi(s)}\le 0 
\label{Gamma1.5.1}
\end{equation}
satisfies, are $(k,s=7)$ where either $7\le k\le 241$ is prime, or
$k=8$, $9$, $16$, $25$, $27$, $32$, $49$; $(k,s=8)$ where $k=8$, $9$, $11$, 
$13$, $17$;
$(k,s=9)$ where $k=9$, $11$, $13$, $17$, $19$; $(k,s=11)$ where $11\le k\le
31$ is prime; $(k,s=13)$ where $13\le k\le 23$ is prime;
$(k=17,s=17)$.

We can take $K_0=911$ in \eqref{case2.6}. Then (here we take
$s_0=6$)
$$
\Delta_1(8)=\ln{\sqrt{2}}- \frac{\ln{7}}{6}-\frac{\ln{911}}{910}
\ge 0.0147667,
$$
and $K_1=38563$ can be taken in \eqref{case2.10}. Checking
\eqref{case2.3} for $6\le s\le k<38563$, we obtain that $6\le s\le
330$ and $6\le s\le k\le 5460$. Moreover, $s\le k\le 330$ for
$20\le s\le 330$. For all these pairs $(k,s)$ satisfying
\eqref{case2.3} which is
\begin{equation}
\frac{\varphi([k,s])}{2\rho(k,s)}\cdot \left(
\ln{\sqrt{2}}-\frac{\ln{\gamma(k)}}{\varphi(k)}-
\frac{\ln{\gamma(s)}}{\varphi(s)}\right) < \ln{\sqrt{4\cdot
14^3}}-\ln{\sin{\frac{\pi}{k}}}- \ln{\sin{\frac{\pi}{s}}}\, ,
\label{Gamma1.5.2}
\end{equation}
we obtain
\begin{equation}
[\bff_{k,s}:\bq]=\frac{\varphi([k,s])}{2\rho(k,s)}\le 909
\label{Gamma1.5.3}
\end{equation}
where $909$ is achieved for $(k,s)=(607,7)$. Moreover, for all
these non-exceptio\-nal pairs $(k,s)$ we obtain the bound
\eqref{case2.11} which is
\begin{equation}
N_0=[\bk:\bff_{k,s}]\le \left[\frac
{\ln{\sqrt{2}}-\ln{\sin{\frac{\pi}{k}}}- \ln{\sin{\frac{\pi}{s}}}}
{\frac{\varphi([k,s])}{2\rho(k,s)}\cdot \left( \ln{\sqrt{4\cdot
14^3}}-\frac{\ln{\gamma(k)}}{\varphi(k)}-
\frac{\ln{\gamma(s)}}{\varphi(s)}\right)}\right]
\label{Gamma1.5.4}
\end{equation}
and finally we obtain the bound \eqref{case2.12} which is
\begin{equation}
N=[\bk:\bq]\le \left[\frac
{\ln{\sqrt{2}}-\ln{\sin{\frac{\pi}{k}}}- \ln{\sin{\frac{\pi}{s}}}}
{\frac{\varphi([k,s])}{2\rho(k,s)}\cdot \left( \ln{\sqrt{4\cdot
14^3}}-\frac{\ln{\gamma(k)}}{\varphi(k)}-
\frac{\ln{\gamma(s)}}{\varphi(s)}\right)}\right] \cdot
\frac{\varphi([k,s])}{2\rho(k,s)}\ . \label{Gamma1.5.5}
\end{equation}
If either a pair $(k,s)$ is exceptional, or the right hand side of
\eqref{Gamma1.5.5} is more than $909$ (these are possible only for
pairs $(k,s)$ with $6\le s\le 17$ and $s\le k\le 421 $), we also
apply to the pair $(k,s)$ the method A of the Case 2 to improve
the poor bound \eqref{Gamma1.5.4} of $N_0=[\bk:\bff_{k,s}]$ for
non-exceptional $(k,s)$. We can apply this method to any pair
$(k,s)$ with $k\ge s\ge 6$ since \eqref{case2.14} is valid if
$a_1+a_2=10$. We obtain that $[\bk:\bq]\le 909$ for all $k\ge s\ge
6$.

Let us assume that $s=3,4$ or $5$ is exceptional. It means that
either $m_1=3,4,5$ or $m_3=3,4,5$ for $\Gamma^{(5)}_1(14)$. For
example, let $m_1=3,4,5$. Let us consider the V-arithmetic graph
defined by $e$, $\delta_1$ and $\delta_4$. We denote
$\alpha=a_{14}^2$ where the algebraic integer
$a_{14}=\delta_1\cdot \delta_4$. The determinant of the Gram
matrix of $e$, $\delta_1$ and $\delta_4$ is equal to
$2\alpha-8\sin^2{(\pi/m_1)}$. It follows that
$0<\sigma(\alpha)<4\sin^2{(\pi/m_1)}$ for
$\sigma\not=\sigma^{(+)}$, and $4<\sigma^{(+)}(\alpha)<14^2$. Then
$\bk=\bq(\alpha)$ and $\bff_{m_1}\subset \bk$. Thus, we can apply
the method A of Case 1 in Sec. \ref{subsec:genresults} to $a_1=0$,
$a_2=4$, $b_1=4$, $b_2=14^2$ and $l:=m_1$ where $m_1=3,4,5$. We
obtain that $[\bk:\bq]\le 76$ for $m_1=3$, $[\bk:\bq]\le 31$ for
$m_1=4$, and $[\bk:\bq]\le 24$ for $m_1=5$.

Thus, finally, $[\bk:\bq]\le 909$ for all graphs
$\Gamma^{(5)}_1(14)$.


\subsection{V-arithmetic 3-graphs $\Gamma_4^{(6)}(14)$
and their ground fields}
\label{subsec:3graph46}

Here we consider V-arithmetic 3-dimensional graphs
$\Gamma_4^{(6)}(14)$ and their fields.

\begin{figure}
\begin{center}
\includegraphics[width=11cm]{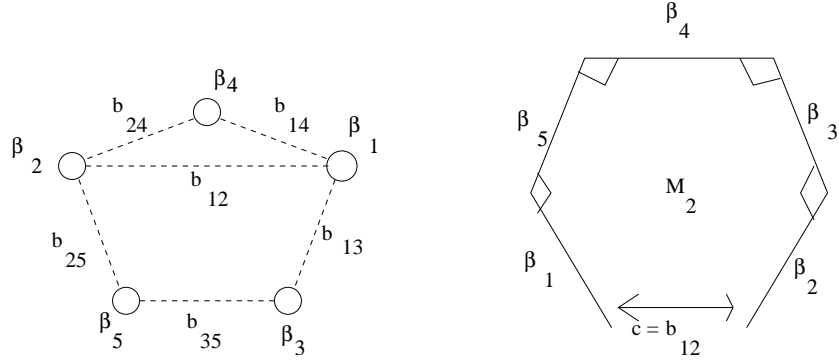}
\end{center}
\caption{The graph of $\M_2$ for $\Gamma_4^{(6)}$}
\label{plgraph46}
\end{figure}

First, let us consider the corresponding plane graph defined by
$\beta_1=\delta_1$, $\beta_2=\delta_2$,
$\beta_3=\delta_3$, $\beta_4=\widetilde{\delta}_4$, $\beta_5=\delta_5$
which give $P(\M_2)$.
We denote $b_{ij}=\beta_i\cdot \beta_j$ when it is not $0$. We also
denote $c=b_{12}$. One can consider it as kind of angle between the
corresponding lines.
This graph is given in Figure \ref{plgraph46}.

Considering determinants of Gram matrices of subsets of $\beta_1,\dots ,
\beta_5$, we obtain all equations of $\M_2$:
\begin{equation}
\left\{
\begin{array}{l}
4b_{13}^2=(4-b_{14}^2)(4-b_{35}^2)\\
4b_{25}^2=(4-b_{24}^2)(4-b_{35}^2)\\
4b_{14}^2+4c^2+4cb_{14}b_{24}=(4-b_{13}^2)(4-b_{24}^2)\\
4b_{24}^2+4c^2+4cb_{14}b_{24}=(4-b_{14}^2)(4-b_{25}^2)\\
b_{35}^2(4-c^2)+2cb_{13}b_{25}b_{35}+4c^2=(4-b_{13}^2)(4-b_{25}^2).
\end{array}
\right.
\label{equationgam46}
\end{equation}

For $\sigma\not=\sigma^{(+)}$ they also satisfy inequalities:
$b_{ij}^2<4$ for all $b_{ij}$, $c^2<4$. By direct calculation of
minimum and maximum of $b_{13}b_{14}b_{24}b_{25}b_{35}$ for
$0\le b_{ij}^2\le 4$ and $0\le c^2\le 4$ satisfying equations
\eqref{equationgam46}, we obtain that
\begin{equation}
-3.1<\sigma(b_{13}b_{14}b_{24}b_{25}b_{35})<3.1\ .
\label{bineqgam46}
\end{equation}
Here minimum $-3.07\dots$ is achieved for $c=-0.39\dots$,
$b_{14}=b_{24}=\pm 1.4\dots$, $b_{13}=\pm 1.166\dots$, $b_{25}=\pm
1.166\dots$, $b_{35}=\pm 1.1549\dots$ and $c=0.39\dots$,
$b_{14}=-b_{24}=\pm 1.4\dots$, $b_{13}=\pm 1.166\dots$,
$b_{25}=\pm 1.166\dots$, $b_{35}=\pm 1.1549\dots$. Here maximum
$3.07\dots$ is achieved for $c=-1.569\dots$, $b_{14}=b_{24}=\pm
1.4\dots$, $b_{13}=\pm 1.166\dots$, $b_{25}=\pm 1.166\dots$,
$b_{35}=\pm 1.1549\dots$ and $c=1.569\dots$, $b_{14}=-b_{24}=\pm
1.4\dots$, $b_{13}=\pm 1.166\dots$, $b_{25}=\pm 1.166\dots$,
$b_{35}=\pm 1.1549\dots$ .

From expressions of $\beta_i$ using $\delta_i$ and $e$, we get
\begin{equation}
b_{13}b_{14}b_{24}b_{25}b_{35}=\frac
{a_{13}a_{14}a_{24}a_{25}a_{35}}
{\sin^2{\frac{\pi}{m}}}.
\label{bagam46}
\end{equation}
We consider the algebraic integer $\alpha\in \bk$ which is
\begin{equation}
\alpha=a_{13}a_{14}a_{24}a_{25}a_{35}.
\label{alpha46}
\end{equation}
From \eqref{bineqgam46}, we obtain
\begin{equation}
-3.1\cdot \sigma(\sin^2{\frac{\pi}{m}})<\sigma(\alpha)<
3.1\cdot \sigma(\sin^2{\frac{\pi}{m}}).
\label{ineqalpha}
\end{equation}
For the geometric embedding $\sigma^{(+)}$ we have that
\begin{equation}
2^5<\sigma^{(+)}(\alpha)<14^5.
\label{ineqalpha+}
\end{equation}
It follows that $\bk=\bq(\alpha)$.

We can apply the methods A and B of Case 1 in Sec.
\ref{subsec:genresults} to $a_1=3.1$, $a_2=3.1$ (then $a=3.1$),
$b_1=2^5$, $b_2=14^5$ (then $b=14^5$), and $l:= m$.

At first, we apply the Method B. All exceptional $l\ge 3$ that is
when \eqref{case1.2} which is
\begin{equation}
\ln{\frac{2}{\sqrt{3.1}}-\frac{\ln{\gamma(l)}}{\varphi(l)}}\le 0
\label{Gamma2.6.1}
\end{equation}
satisfies are $l=3,4,5,7,8,9,11,13,17,19,23$.

We can take $L_0=2053$ in \eqref{case1.6}. Then
$$
\Delta_1(3.1)= \ln(\frac{2}{\sqrt{3.1}})-\frac{\ln{2053}}{2052}>
0.1237\,,
$$
and $L_1=2125$ can be taken in \eqref{case1.9}. Checking
\eqref{case1.3} for $3\le  l< 2125$, we obtain that $3\le l\le
510$. For all these $l$ such that \eqref{case1.3} which is
\begin{equation}
\frac{\varphi(l)}{2}\cdot \left(\ln{\frac{2}{\sqrt{3.1}}-
\frac{\ln{\gamma(l)}}{\varphi(l)}}\right)
<\ln{\sqrt{\frac{14^5}{3.1}}}-\ln{\sin{\frac{\pi}{l}}}
\label{Gamma2.6.2}
\end{equation}
satisfies, we obtain
\begin{equation}
[\bff_{l}:\bq]=\frac{\varphi(l)}{2}\le 99 \label{Gamma2.6.3}
\end{equation}
where $99$ is achieved for $l=199$. Moreover, for all these
non-exceptional $l$ we obtain the bound \eqref{case1.10} which is

\begin{equation}
N_0=[\bk:\bff_l]\le \left[
\frac{\ln{\sqrt{\frac{14^5}{3.1}}}-\ln{\sin{\frac{\pi}{l}}}}
{\frac{\varphi(l)}{2}\cdot \left(\ln{\frac{2}{\sqrt{3.1}}-
\frac{\ln{\gamma(l)}}{\varphi(l)}}\right)          } \right],
\label{Gamma2.6.4}
\end{equation}
and finally we obtain the bound \eqref{case1.11} which is
\begin{equation}
N=[\bk:\bq]\le \left[
\frac{\ln{\sqrt{\frac{14^5}{3.1}}}-\ln{\sin{\frac{\pi}{l}}}}
{\frac{\varphi(l)}{2}\cdot \left(\ln{\frac{2}{\sqrt{3.1}}-
\frac{\ln{\gamma(l)}}{\varphi(l)}}\right)          } \right]\cdot
\frac{\varphi(l)}{2}\ . \label{Gamma2.6.5}
\end{equation}
If either $l$ is exceptional, or the right hand side of
\eqref{Gamma2.6.5} is more than $99$ (this is possible for $3\le
l\le 113$ only), we also apply to $l$ the method A of the Case 1
to improve the poor bound \eqref{Gamma2.6.4} for
$N_0=[\bk:\bff_{l}]$ for non-exceptional $l$. We can apply this
method to any $l\ge 4$ since \eqref{case1.12a} is valid for all
$l\ge 4$ if $a_1+a_2=6.2$. For $l\ge 6$, this method gives what we
want: $[\bk:\bq]\le 99$. For $l=4$, it only gives $[\bk:\bq]\le
120$; for $l=5$, it only gives $[\bk:\bq]\le 172$.

If $l=3,4,5$, equivalently $m=3,4,5$, considering the subgraph of
$e, \delta_4, \delta_2$, exactly the same consideration as for the
graph $\Gamma^{(5)}_1(14)$ above for $m_1=3,4,5$, give that
$[\bk:\bq]\le 76$ for $m=3$, $[\bk:\bq]\le 31$ for $m=4$, and 
$[\bk:\bq]\le 24$ for $m=5$.

Thus, $[\bk:\bq]\le 99$ for all graphs $\Gamma_4^{(6)}(14)$.

\medskip

This finishes the proof of Theorem \ref{th3graphs}.

\section{Appendix: Some results about\\ cyclotomic fields} \label{appendix1}
This is exactly the same as in \cite{Nik7}.  We repeat it for
readers convenience.

Here we give some results about cyclotomic fields which we used. All of
them follow from standard results. For example, see the book \cite{CF}.

We consider the cyclotomic field $\bq\left(\sqrt[l]{1}\right)$ and
its totally real subfield $\bff_l=\bq\left(\cos{(2\pi/l)}\right)$.
We have $[\bq\left(\sqrt[l]{1}\right):\bq]=\varphi(l)$ where
$\varphi(l)$ is the Euler function. We have
$\bff_l=\bq\left(\sqrt[l]{1}\right)=\bq$ for $l=1,2$, and
$[\bff_l:\bq]=\varphi(l)/2$ for $l\ge 3$. It is known (e.g., see
\cite{CF}) that the discriminant of the field
$\bq\left(\sqrt[l]{1}\right)$ is equal to (where $p$ is prime)
\begin{equation}
|{\rm discr}\,\bq(\sqrt[l]{1})|=\frac{l^{\varphi(l)}}{\prod_{p|l}
{p^{\varphi(l)/(p-1)}}}\ . \label{discrql}
\end{equation}

Let $\zeta_l=\exp{(2\pi i/l)}$ be a primitive $l$-th root of $1$.
The element $\zeta_l$ generates the ring of integers of
$\bq\left(\sqrt[l]{1}\right)$. Further we assume that $l\ge 3$.
The equation of $\zeta_l$ over $\bff_l$ is
$g(x)=(x-\zeta_l)(x-\zeta_l^{-1})=x^2-(\zeta_l+\zeta_l^{-1})x+1=0$.
We have
$g^\prime(\zeta_l)=2\zeta_l-(\zeta_l+\zeta_l^{-1})=\zeta_l-\zeta_l^{-1}$.
Thus,
$$
N_{\bq\left(\sqrt[l]{1}\right)/\bff_l}(g^\prime(\zeta_l))=
(\zeta_l-\zeta_l^{-1})(\zeta_l^{-1}-\zeta_l)=4\sin^2{(2\pi/l)}
$$
which gives the discriminant $\delta
\left(\bq\left(\sqrt[l]{1}\right)/\bff_l\right)=4\sin^2{(2\pi/l)}$.
It follows
$$
|\delta\left(\bq\left(\sqrt[l]{1}\right)/\bq\right)|=
|\delta\left(\bff_l/\bq\right)^2N_{\bff_l/\bq}(4\sin^2{(2\pi/l)})|.
$$
We have
\begin{equation}
N_{\bff_l/\bq }(4\sin^2{(\pi/l)})=\gamma(l)=\left\{
\begin{array}{cl}
p &\ {\rm if}\  l=p^t>2\ {\rm where}\ p\ {\rm is}\ {\rm prime ,} \\
1 &\ {\rm otherwise.}
\end{array}\
\label{normsin}
 \right .
\end{equation}

If $l$ is odd, then $4\sin^2{(\pi/l)}$ and $4\sin^2{(2\pi/l)}$ are
conjugate, and their norms are equal. Thus,
$$
N_{\bff_l/\bq}(4\sin^2{(2\pi/l)})=\gamma(l),\ \text{if $l\ge 3$ is
odd.}
$$

If $l$ is even and $l_1=l/2$, then
$4\sin^2{(2\pi/l)}=4\sin^2{(\pi/l_1)}$. If $l_1$ is odd, then
$\bff_{l_1}=\bff_{l}$, and we get
$$
N_{\bff_l/\bq}(4\sin^2{(2\pi/l)})=\gamma(l/2)\ \text{if $l\ge 6$
is even, but $l/2$ is odd.}
$$
If $l_1\ge 4$ is even, then $[\bff_l:\bff_{l_1}]=2$, and we get
$$
N_{\bff_l/\bq}(4\sin^2{(2\pi/l)})=\gamma(l/2)^2\ \text{if $l\ge 8$
and $l/2$ is even.}
$$
At last,
$$
N_{\bff_4/\bq}(4\sin^2{(2\pi/4)})=4
$$
if $l=4$.

Thus, finally we get for $l\ge 3$:
\begin{equation}
N_{\bff_l/\bq }(4\sin^2{(2\pi/l)})=\widetilde{\gamma}(l)=\left\{
\begin{array}{cl}
\gamma(l) &\ {\rm if}\ l\ge 3\ {\rm is\ odd,} \\
\gamma(l/2) &\ {\rm if}\ l/2 \ge 3\ {\rm is\ odd,}\\
\gamma(l/2)^2 &\ {\rm if}\ l/2\ge 4\ {\rm is\ even,}\\
4 &\ {\rm if}\ l=4. \\
\end{array}\
\label{normsin2}
 \right .
\end{equation}
Moreover, we obtain the formula for the discriminant:
\begin{equation}
|{\rm discr}\,\bff_l|= \left(|{\rm discr}\,\bq
(\sqrt[l]{1})|\left/\right.\widetilde{\gamma}(l)\right)^{1/2}\
\text{for}\ l\ge 3
\label{discrFl}
\end{equation}
where $|{\rm discr}\,\bq (\sqrt[l]{1})|$ is given by
\eqref{discrql}, and $\widetilde{\gamma}(l)$ is given by
\eqref{normsin2}.

We denote
$\bff_{k,s}=\bq\left(\cos{(2\pi/k)},\,\cos{(2\pi/s)}\right)$.
Further we assume that $k,s\ge 3$. Let $m=[k,s]$ be the least
common multiple of $k$ and $s$. Then $\bff_{k,s}\subset
\bff_m\subset \bq(\sqrt[m]{1})$. We have
$\text{Gal}\,\left(\bq(\sqrt[m]{1})/\bq\right)=(\bz/m\bz)^\ast$
where $\alpha\in (\bz/m\bz)^\ast$ acts on each $m$-th root $\zeta$
of $1$ by the formula $\zeta\mapsto \zeta^\alpha$. Obviously,
$\bff_{k,s}$ is the fixed field of the subgroup $G$ of the Galois
group $(\bz/m\bz)^\ast$ which consists of all $\alpha \in
(\bz/m\bz)$ such that $\alpha\equiv \pm 1\mod k$ and $\alpha\equiv
\pm 1\mod s$. The $G$ includes the subgroup of order two of
$\alpha\equiv \pm 1\mod m$. If $\alpha\equiv 1\mod k$, then
$\alpha\equiv 1+kt\mod m$ where $t\in \bz$. If $1+kt\equiv -1\mod
s$, then the equation $kt+sr=2$ has an integer solution $(t,r)$
which is equivalent to $(k,s)|2$. Thus, $G$ has the order $4$ if
and only if  $(k,s)|2$. Otherwise, $G$ has the order $2$. We set
\begin{equation}
\rho(k,s)=\left\{
\begin{array}{cl}
2 &\ {\rm if}\ (k,s)|2, \\
1 &  {\rm otherwise,}\\
\end{array}\
\label{defrho}
 \right.
\end{equation}
and we obtain
\begin{equation}
[\bff_{k,s}:\bq]=\frac{\varphi(m)}{2\rho(k,s)}. \label{degreekl}
\end{equation}
Moreover, we get
$$
\bff_{k,s}=\bff_m\ \text{if}\ (k,s)\not|\, 2.
$$
It follows,
\begin{equation}
|\text{discr}\,\bff_{k,s}|=|\text{discr}\,\bff_m|\ \text{if}\
(k,s)\not|\,2, \label{discrFks1}
\end{equation}
where $m=[k,s]$, and $|\text{discr}\,\bff_m|$ is given by
\eqref{discrFl}.

Assume that $(k,s)|2$. If $(k,s)=1$, then the fields
$\bq(\sqrt[k]{1})$ and $\bq(\sqrt[s]{1})$ are linearly disjoint
and their discriminants are coprime. Then their subfields $\bff_k$
and $\bff_s$ are linearly disjoint, and their discriminants are
coprime, and we obtain
$$
|\text{discr}\,\bff_{k,s}|=|\text{discr}\,\bff_k|^{(\varphi(s)/2)}
|\text{discr}\,\bff_s|^{(\varphi(k)/2)}\ \text{if}\ (k,s)=1\
\text{and}\ k,s\ge 3.
$$

Assume that $(k,s)=2$. Then one of $k/2$ or $s/2$ is odd. Assume,
$k_1=k/2$ is odd. Then $\bff_k=\bff_{k_1}$ and
$\bff_{k,s}=\bff_{k_1,s}$ where $(k_1,s)=1$. Thus, we obtain the
previous case which gives exactly the same formula. We finally
obtain
\begin{equation}
|\text{discr}\,\bff_{k,s}|=|\text{discr}\,\bff_k|^{(\varphi(s)/2)}
|\text{discr}\,\bff_s|^{(\varphi(k)/2)}\ \text{if}\ (k,s)|2\
\text{and}\ k,s\ge 3 \label{discrFks2}
\end{equation}
where the discriminants $|\text{discr}\,\bff_k|$ and
$|\text{discr}\,\bff_s|$ are given by \eqref{discrFl}.


V.V. Nikulin \par Deptm. of Pure Mathem. The University of
Liverpool, Liverpool\par L69 3BX, UK; \vskip1pt Steklov
Mathematical Institute,\par ul. Gubkina 8, Moscow 117966, GSP-1,
Russia

vnikulin@liv.ac.uk \ \ vvnikulin@list.ru

\end{document}